\documentclass[fleqn]{mat01}
\usepackage{times,mathtimy,amssymb,latexsym}
\begin{document}

\def\d{\mbox{\rm d}}

\newtheorem{theore}{Theorem}
\renewcommand\thetheore{\arabic{section}.\arabic{theore}}
\newtheorem{theor}{\bf Theorem}
\newtheorem{lem}{Lemma}
\newtheorem{exampl}{Example}

\def\remar{\trivlist \item[\hskip \labelsep{\it Remarks.}]}
\def\rema{\trivlist \item[\hskip \labelsep{\it Remark.}]}

\renewcommand\theequation{\arabic{equation}}

\setcounter{page}{383}
\firstpage{383}

\title{An algebra of absolutely continuous functions\\ and its multipliers}

\markboth{Savita Bhatnagar}{Absolutely continuous functions and its
multipliers}

\author{SAVITA BHATNAGAR}

\address{Department of Mathematics, Panjab University, Chandigarh~160~014, India\\
\noindent E-mail: bhsavita@pu.ac.in}

\volume{115}

\mon{November}

\parts{4}

\pubyear{2005}

\Date{MS received 28 February 2005; revised 9 September 2005}

\begin{abstract}
The aim of this paper is to study the algebra $AC_{p}$ of absolutely
continuous functions $f$ on [0,1] satisfying $f(0) = 0, f'\in
L^{p} [0, 1]$ and the multipliers of $AC_{p}$.
\end{abstract}

\keyword{Absolutely continuous function; Banach algebra; multiplier.}

\maketitle

\section{Introduction}

Let $I = [0, 1]$ be the compact topological semigroup with max
multiplication and usual topology. $C(I), L^{p}(I), 1\leq
p\leq\infty$ are the associated Banach algebras. Larsen \cite{5}
obtained multipliers for the Banach algebra $L^{1}(I)$. Baker, Pym
and Vasudeva \cite{1} obtain characterizations of multipliers from
$L^{p}(I)$ to $L^{r}(I), 1\leq r, p\leq\infty$. Bhatnagar and
Vasudeva \cite{2} characterize $\hbox{Hom}_{C(I)}(L^{r}(I),
L^{p}(I)) $ and their pre-duals for $1\leq r\leq p.$ Bhatnagar
\cite{3} studied $\hbox{Hom}_{C(I)}(L^{r}(I), L^{p}(I))$ and their
pre-duals for $r>p$. The study of pre-duals of the multipliers in
\cite{2} and \cite{3} involved a deep understanding of the
interpolation theory and lengthy calculations. The multipliers can
be obtained more attractively by abstract arguments. It turns out
that results obtained via abstract arguments \hbox{compare} with
those obtained in \cite{2} and \cite{3}. Our results also include
the results obtained by Larsen \cite{5}. $L^{p}(I)$ is replaced by
the Gelfand transform algebra $AC_{p}, p\geq 1$ of absolutely
continuous functions $f$ on [0,1] with $f(0) = 0$ and $f{'}\in
L^{p}(I).$ With the norm defined as $|||f|||=\|f{'}\|_{p}, f\in
AC_{p}, AC_{p}$ constitute subalgebras of $C(I)$ and have an
approximate identity. The purpose of this note is to study the
multipliers from $AC_{r}$ to $AC_{p}, 1\leq r, p\leq\infty$. For
$r = p$, a complete description of multipliers is obtained. For $r
< p$, the multiplier algebra consists of $\{0\}$ alone. In case of
$r > p$ we provide a set of necessary and another set of
sufficient conditions for a function in $C(0,1]$ to be a
multiplier. An example of a function which satisfies the necessary
conditions but does not satisfy the sufficient conditions and
fails to be a multiplier is also provided.

\section{The Banach algebra $\pmb{AC_{p}}$}

Let $I = [0, 1]$ with the usual interval topology be the compact metric
space and $C(I)$ be the set of all continuous complex valued functions
on $I$. For $1\leq p\leq\infty$, let
\begin{align*}
AC_{p} &= \{f\in C(I)\hbox{\rm :}\ f\ \hbox{is absolutely continuous on} \
I, f(0) = 0\\[.2pc]
&\quad\, \hbox{and} \ f{'}\in L^{p}(I)\}.
\end{align*}

Define $|||f|||=\|f{'}\|_{p}, f\in AC_{p}.$ For $1\leq p\leq\infty$,
let $p{'}$ denote the conjugate index to $p$, that is, $\frac{1}
{p}+\frac{1}{p'}=1$.

The following inequality will be frequently used in the sequel. For $g
\in AC_{p}, 1 < p < \infty,$
\begin{align}
&|g(t)| = \left|\int^{t}_{0} g'(s) \d s \right|\leq
\int^{t}_{0}|g'(s)| \d s \leq \left(\int^{t}_{0} |g'(s)|^{p} \d s
\right)^{1/p} t^{1/p'} \leq \|g'\|_{p},\nonumber\\[.2pc]
&\hbox{for} \ g\in AC_{1}, |g(t)|\leq \|g'\|_{1}, \ \hbox{and for}
\ g\in AC_{\infty}, |g(t)|\leq \|g'\|_{\infty}.t,
\end{align}
using Holder's inequality.

\begin{lem}
Let $f$ be continuous on $[0, 1]${\rm ,} absolutely continuous on $[\alpha, 1]$
for any $\alpha$ with $0 < \alpha \leq 1$ with $f'\in L^{1} [0, 1]$. Then
$f$ is absolutely continuous on $[0,1]$.
\end{lem}

The proof is elementary.

\setcounter{theor}{1}
\begin{theor}[\!]
$AC_{p}$ is a Banach subalgebra of $C(I)$. It has an approximate
identity for $1\leq p < \infty$. $AC_{\infty}$ has no approximate
identity. The maximal ideal space $\triangle(AC_{p})$ of $AC_{p},
1\leq p\leq \infty$ is homeomorphic to {\rm (0,1]}. Moreover{\rm,}
$AC_{r}\subseteq AC_{p}$ for $r > p$.
\end{theor}

\begin{proof}
Clearly $AC_{p}$ is a Banach space. For $f, g \in AC_{p}$,
\begin{equation*}
(fg)' = f'g + fg'\in L^{p}(I)
\end{equation*}
as $f', g'\in L^{p}(I)$ and $f, g$ are bounded. Moreover,
\begin{equation*}
|||fg||| = \|(fg)'\|_{p} \leq \|g\|_{\infty} \|f'\|_{p} + \|f\|_
{\infty} \|g'\|_{p}.
\end{equation*}

Since $g \in AC_{p}, |g(t)| \leq \|g'\|_{p} = |||g|||$, using (1),
and so $\|g\|_{\infty} \leq |||g|||.$

Similarly $\|f\|_{\infty} \leq
|||f|||.$ Consequently, $|||fg|||\leq 2|||f|||~|||g|||.$ Thus $AC_{p}$
is a Banach algebra.

Define $e_{\alpha}(t) = \min\{t/\alpha, 1\}, t\in I.$ Then
$e'_{\alpha}(t) = \frac{1}{\alpha}\chi_{[ 0,\alpha)}(t)$. We show that
$\{e_{\alpha}\}$ forms an approximate identity for $AC_{p}, 1\leq
p<\infty$ as $\alpha \rightarrow 0+.$

Now
\begin{align*}
|||e_{\alpha} g-g||| &=\|(e_{\alpha}-1)g'+e'_{\alpha}g\|_{p}\\[.2pc]
 &\leq\|(e_{\alpha}-1)g'\|_{p} + \|e'_{\alpha}g\|_{p}.\\[.2pc]
\|e'_{\alpha}g\|_{p} &= \left(\int^{\alpha}_{0} \left|\frac{1}{\alpha}g(t)\right|^{p} \d t\right)^{1/p}\\[.2pc]
&\leq\frac{1}{\alpha}\left(\int^{\alpha}_{0} \left(\int^{t}_{0}|g'(s)|^{p} \d s\right) t^{p/p'} \d t\right)^{1/p},\quad \hbox{using (1)}.
\end{align*}

So
\begin{align*}
\|e'_{\alpha}g\|_{p} &\leq \frac{1}{\alpha}\|g'\chi_{[0,\alpha]}\|_{p}
\left(\int^{\alpha}_{0} t^{p/p'} \d t\right)^{1/p}\\[.2pc]
&= \frac{1}{\alpha} \|g'\chi_{[0,\alpha]}\|_{p} \left(\frac{\alpha^{p}}{p}\right)^{1/p}\\[.2pc]
&= p^{-1/p} \|g'\chi_{[0,\alpha]}\|_{p}.
\end{align*}

Also\vspace{-1pc}
\begin{align*}
\|(e_{\alpha}-1)g'\|_{p} &= \left(\int^{\alpha}_{0} \left|\frac{t}{\alpha}-1\right|^{p}
  |g'(t)|^{p} \d t\right)^{1/p}\\
&\leq \|g'\chi_{[0,\alpha]}\|_{p} \quad \hbox{as} \ \
\left|\frac{t}{\alpha}-1 \right|\leq 1 \quad \hbox{for} \ \ t\in
[0,\alpha].
\end{align*}

Thus $|||e_{\alpha} g-g||| \leq (p^{-1/p} + 1)\|g'\chi_{[0,\alpha]}\|_{p} \rightarrow 0$ as $\alpha \rightarrow 0 +$
because $g'\in L^{p}(I), 1\leq p < \infty.$

Next, if $g \in AC_{\infty}$ and $\|g'\|_{\infty} = K$ then
$|g(t)|=|\int^{t}_{0}g'(s) \d s| \leq Kt$ and $|(tg)'|=|tg'+g| \leq 2Kt$
for all $t\in I$. So
\begin{equation*}
|||tg-t||| = \|tg'+g-1\|_{\infty} \geq \lim_{t\rightarrow 0} |1-2Kt|=1.
\end{equation*}

So $AC_{\infty}$ has no approximate identity. We now find the maximal
ideal space $\triangle(AC_{p})$ of $AC_{p}$.

Clearly, $AC_{p}$ separates strongly the points of [0,1] and is
self-adjoint. If $f\in AC_{p}$ and $a = \hbox{inf}_{t}|1 - f(t)| > 0$, take
$g = \frac{f}{f-1}$. Then $g' = -\frac{f'}{(f-1)^{2}}$, so that $|g'|\leq
\frac{|f'|}{a^{2}} \in L^{p}(I)$, i.e., $g\in AC_{p}$ and
$g\!\circ\!f = g+f-gf=0.$ Thus $f$ is quasiregular in $AC_{p}$. It follows using
Corollary~3.2.8 of Rickart \cite{7} that the maximal ideal space
$\triangle(AC_{p})$ of $AC_{p}$ is homeomorphic to $\triangle (C(I))=I$
under the natural embedding. Since for $f\in AC_{p}, \hat f(0)= f(0)=0$
we get that $\triangle(AC_{p})$ of $AC_{p}$ is homeomorphic to (0,1].

Finally, for $r > p, f\in AC_{r}$ we have $f'\in L^{r}(I)\subset
L^{p}(I)$, so that $f\in AC_{p}$. The inclusion $AC_{r} \subset AC_{p}$
is indeed proper as $L^{r}(I)$ is a proper subset of $L^{p}(I)$ for
$r>p$. This completes the proof. \hfill $\Box$
\end{proof}

Note that the approximate identity $\{e_{\alpha}\}$ is bounded if $p=1$
and is unbounded if $1 < p < \infty$. Also with our methods the case
$p = \infty$ has been solved completely whereas in \cite{1} the maximal
ideal space of $L^{\infty}(I)$ could not be calculated.\vspace{-.5pc}

\section{The multiplier space}

A mapping $T$ on a commutative Banach algebra $A$ to itself is called a
multiplier if $T(xy) = xT(y) = T(x)y, x, y \in A$. If $A$ is semisimple and
$T\hbox{\rm :}\ A\rightarrow A$ is a multiplier then there exists a unique bounded
continuous function $m$ on $\triangle (A)$ such that $\hat Tx= m \hat x$
for all $x\in A$ and $\|m\|_{\infty}\leq \|T\|$ (p.~19 of \cite{4}). Since
for $g\in AC_{p}, \hat g(t) = g(t), t \in (0,1], AC_{p}$ is a semisimple
Banach algebra, a multiplier $T\hbox{\rm :}\ AC_{p}\rightarrow AC_{p}$ is a map
satisfying $Tg = mg, g\in AC_{p}$ for some continuous bounded function
$m$ on (0,1]. The following theorem gives necessary and sufficient
conditions for $m\in C_{b}(0,1]$ to be a multiplier of $AC_{p}$.

\begin{theor}[\!]
A map $T\hbox{\rm :}\ AC_{p}\rightarrow AC_{p}, 1\leq p\leq \infty$ is a
multiplier iff there exists an $m \in C_{b}(0,1]$ such that for each
$\epsilon >0, m$ is absolutely continuous on $[\epsilon,1], m'\in
L^{p}[\epsilon,1]$ and $\|m'\chi_{[\epsilon,1]}\|_{p} = O(\epsilon^{-1/p'}).$ {\rm (}Treat
$\epsilon^{-1/p'}=1$ for $p=1$.{\rm )}
\end{theor}

\begin{proof}
Suppose $T$ is a multiplier of $AC_{p}$. Then there exists $m\in
C_{b}(0,1]$ such that $Tg = mg$. Let the norm of the
multiplication operator be $N$. Then $|||Tg|||=|||m g|||\leq
N|||g|||, g \in AC_{p}.$

As $m e_{\alpha}\in AC_{p}$ and $e_{\alpha}=1$ on $[\alpha,1]$, we
get that $m$ is absolutely continuous on $[\alpha,1], m'\in
L^{p}[\alpha,1]$ and $|||m e_{\alpha}|||\leq N |||e_{\alpha}|||=
N\|e'_{\alpha}\|_{p}= N\alpha^{-1/p'}\!.$

Also for $1 \leq p < \infty$,
\begin{align*}
|||m e_{\alpha}|||^{p} &= \|m' e_{\alpha} + m e'_{\alpha}\|^{p}_{p}\\[.2pc]
&= \int^{\alpha}_{0}|m' e_{\alpha} + m e'_{\alpha}|^{p}(t) \d t + \int_{\alpha}^{1}|m'
e_{\alpha} + m e'_{\alpha}|^{p}(t) \d t\\[.2pc]
&\geq \int_{\alpha}^{1}|m'|^{p}(t) \d t \quad \hbox{as} \
e_{\alpha}=1 \ \hbox{and} \  e'_{\alpha}=0 \ \hbox{on} \ [\alpha,1].
\end{align*}

So $(\int_{\alpha}^{1}|m'|^{p}(t) \d t)^{1/p}\leq N \alpha^{-1/p'}$ \ or
$\|m' \cdot \chi _{[\alpha,1]}\|_{p}=O(\alpha^{-1/p'}).$ For $p=\infty$,
\begin{align*}
|||me_{\alpha}||| &= \|m' e_{\alpha} + me'_{\alpha}\|_{\infty}\\[.2pc]
&\geq\|(m' e_{\alpha} + me'_{\alpha}) \cdot \chi _{[\alpha,1]}\|_{\infty}\\
&= \|m' \cdot \chi _{[\alpha,1]}\|_{\infty}
\end{align*}
and
\begin{equation*}
|||m e_{\alpha}|||\leq N |||e_{\alpha}|||= N\|e'_{\alpha}\|_{\infty}= N/\alpha
\end{equation*}
so that
\begin{equation*}
\|m' \cdot \chi _{[\alpha,1]}\|_{\infty}=O(1/\alpha).
\end{equation*}
Conversely, suppose $m\in C_{b}(0,1]$ satisfies, $m$ is absolutely
continuous on $[\epsilon, 1], m'\in L^{p}[\epsilon, 1],
\|m'\chi_{[\epsilon, 1]}\|_{p} = O (\epsilon^{-1/p'})$ for
$\epsilon >0$ and $g\in AC_{p}$. As $g(0) = 0$ we have that $mg$
has a continuous extension to $[0, 1]$ by assigning $mg(0) = 0$.
We first show that $(m g)'\in L^{p}(I)$. Now $(mg)' = mg' + m'g$.
Since $m\in C_{b}(0, 1]$ and $g'\in L^{p}(I)$, we get that $m
g'\in L^{p}(I)$. It remains to prove that $m' g\in L^{p}(I).$

For $p = 1, p'=\infty$ and $\|m'\chi_{[\epsilon,1]}\|_{1}=O(1)$ implies
that $\lim_{\epsilon \rightarrow 0}\int^{1}_{\epsilon}|m'(t)| \d t$ exists,
so $m'\in L^{1}(I)$. Also $g\in AC_{1}$ is bounded so $m'g \in
L^{1}(I)$ and $\|m'g\|_{1}\leq \|m'\|_{1}~\|g'\|_{1}$, using eq.~(1).

For $1<p<\infty$ and $0<\alpha \leq 1/4$,
\begin{align*}
\int^{2\alpha}_{\alpha}|m'(s)|^{p}|g(s)|^{p} \d s &\leq
\int^{2\alpha}_{\alpha}|m'(s)|^{p} \left(\int ^{s}_{0}|g'(t)| \d t \right)^{p} \d s\\[.1pc]
&\leq \left(\!\int ^{2\alpha}_{0}\!|g'(t)| \d t\!\right)^{p}\!\!\!\cdot A \cdot \alpha^{-p/p'}\!\!,
\hbox{where} \ A \ \hbox{is a constant}\\[.1pc]
&= A \left(\int^{2\alpha}_{0}|g'(t)| \d t \right)^{p} \cdot \alpha^{1-p}\\[.1pc]
&= B \left(\int^{2\alpha}_{0}|g'(t)| \d t \right)^{p} \int^{4\alpha}_{2\alpha} s^{-p}
\d s,\\[.1pc]
&\quad\ \ \hbox{where} \ \ B = \frac{A(p-1)}{2^{1-p}-4^{1-p}} \
\hbox{is a constant,}\\[.1pc]
&\leq B \int ^{4\alpha}_{2\alpha}s^{-p} \left(\int^{s}_{0}|g'(t)| \d t \right)^{p} \d s\\[.1pc]
&= B \int ^{4\alpha}_{2\alpha} \left(\frac{1}{s}\int ^{s}_{0}|g'(t)| \d t \right)^{p} \d s.
\end{align*}

Adding these inequalities for $\alpha=1/4,1/8,1/16, \dots$, we get
\begin{align*}
\int^{1/2}_{0}\!\!\!|m'(s)|^{p}|g(s)|^{p} \d s \leq\!B\!\int
^{1}_{0}\!\!\left(\frac{1}{s}\!\int^{s}_{0}\!|g'(t)| \d
t\!\!\right)^{p}\! \d s \leq B\!\left(\frac{p}{p-1}\|g'\|_{p}\!\!\right)^{p}\!\!<\!\infty
\end{align*}
as $\frac{1}{s}\int^{s}_{0}|g'(t)| \d t \in L^{p}(0,1]$ by Hardy's
inequality \cite{8}. Also, $\int^{1}_{1/2}|m'(s)|^{p}|g(s)|^{p} \d s \leq
\|g'\|^{p}_{p} \cdot A \cdot (1/2)^{-p/p'} <\infty$, using eq.~(1) so that $m'g\in
L^{p}(I)$ and $\|m'g\|_{p} \leq \big(B\big(\frac
{p}{p-1}\big)^{p}+A~(1/2)^{-p/p'}\big)^{1/p} \|g'\|_{p}$.

For $p=\infty$, let $g\in AC_{\infty}$. Then
$|g(t)|=|\int^{t}_{0}g'(s) \d s|\leq t \|g'\|_{\infty}.$

Now $\Big|m'g \chi_{\big(\frac{1}{2^{n}},\frac{1}{2^{n-1}}\big]}
\Big|\leq A 2^{n}\|g'\|_{\infty} \frac{1} {2^{n-1}}=2 A
\|g'\|_{\infty}, n = 1,2,3,\dots$.

Therefore $\|m'g\|_{\infty}\leq 2A\|g'\|_{\infty} <\infty$ or $m'g\in
L^{\infty}(I)$. That $mg$ is absolutely continuous follows from Lemma~1
as $mg$ is absolutely continuous on $[\epsilon,1]$ for all $\epsilon
>0$, continuous at 0 and $(mg)'\in L^{p}(I)\subseteq L^{1}(I)$. Thus $m$
is a multiplier. It is easy to see that the norm of the multiplier $T$
given by $m$ is
\begin{align*}
\|T\| = {\displaystyle \sup_{\|g'\|_{p}=1}}\|m'g+mg'\|_{p}\leq \left\{
\begin{array}{ll}
\|m\|_{\infty}+ \|m'\|_{1}, &p=1\\[.2pc]
\|m\|_{\infty}+\big(B\big(\frac{p}{p-1}\big)^{p} &\\[.2pc]
\quad +A(1/2)^{-p/p'}\big)^{1/p}, &1<p<\infty\\[.2pc]
\|m\|_{\infty}+2A, &p=\infty \end{array}.\right.
\end{align*}

This completes the proof. \hfill $\Box$
\end{proof}

\begin{rema}{\rm
For $p = 1, m\in C_{b}(0,1]$ is a multiplier of
$AC_{1}$ iff $m$ is absolutely continuous on $[\epsilon, 1]$ for
each $\epsilon > 0$ and $m'\in L^{1}(I)$, i.e., iff $m$ is
absolutely continuous on $I$, using Lemma~1 as $m$ being bounded can
be continuously extended to $[0,1]$. Thus Theorem~3 includes the
results of \cite{5}.

Multipliers from $AC_{r}$ to $AC_{p},r\neq p$ are
given by continuous functions $m\in C(0,1]$ which may not be
bounded.}
\end{rema}

\begin{theor}[\!]
If $T\hbox{\rm :}\ AC_{r}\rightarrow AC_{p}$ is a multiplier{\rm ,} $r<p\leq
\infty${\rm ,} then $T = 0$.
\end{theor}

\begin{proof}
Let $T\hbox{\rm :}\ AC_{r}\rightarrow AC_{p}$ be a multiplier, $r<p\leq \infty$.
Then there exists a function $m\in C(0,1]$ such that $Tg=mg$. If $m\neq
0$, then there exists $s_{0}\in (0,1]$ such that $m( s_{0})\neq 0$. By
continuity there is an $\epsilon >0$ such that $|m(t)| > k/2$ for $t\in
(s_{0}-\epsilon, s_{0}+\epsilon)=N$, say. Here $k=|m(s_{0})|$.

For $t\in N$,
\begin{align*}
|(mg)'(t)| &=|(m'g+mg')(t)|\\[.2pc]
&\geq |m(t)g'(t)|-|m'(t)g(t)|\\[.2pc]
&\geq \frac{k}{2}|g'(t)|-|m'(t)g(t)|
\end{align*}
or
\begin{equation*}
\frac {k}{2 }|g'(t)|\leq |(mg)'(t)|+|m'(t)g(t)|. \tag{\hbox{$*$}}
\end{equation*}
For $g = e_{\alpha}, \alpha < s_{0} - \epsilon$ and $t\in N$,
$(mg)'(t)=m'(t)g(t)+m(t)g'(t)=m'(t)$. Thus $m'\in L^{p}(N)$.

From $(*)$, $\frac {k} {2}|g'(t)|\leq |(mg)'(t)|+|m'(t)| \ \
\|g\|_{r}\in L^{p}(N)$ for all $g\in AC_{r}.$ A contradiction. So
$m = 0.$ \hfill $\Box$
\end{proof}

The multipliers from $AC_{r}$ to $AC_{p}, r>p$ are given
by continuous functions on (0,1] which are locally in $AC_{p}$.
The following theorem provides necessary growth conditions on $m'$
for $m$ to be a multiplier from $AC_{r}$ to $AC_{p}, r  > p$.

\begin{theor}[\!]
Let $T\hbox{\rm :}\ AC_{r}\rightarrow AC_{p}, r > p$ be a multiplier given
by $Tg=mg, g\in AC_{r}$. Then for $\epsilon > 0${\rm ,} $m$ is absolutely
continuous on $[\epsilon,1], m'\in L^{p}[\epsilon,1]$ and
$\|m'\chi_{[\epsilon,1]}\|_{p}=O(\epsilon^{-1/r'})$.
\end{theor}

\begin{proof}
(Similar to the proof of Theorem~3.) Let the norm of the
multiplication operator $T$ be $N$. Then $|||Tg|||=|||m g|||\leq
N|||g|||, g \in AC_{r}$ where the norm on the left-hand side is in
$AC_{p}$ and the norm on the right-hand side is in $AC_{r}$.

Since $e_{\alpha} \in AC_{r}$, $\|(m e_{\alpha})'\|_{p}\leq N
\|e'_{\alpha}\|_{r}$ gives (as in Theorem~3)
$\|m'\chi_{[\alpha,1]}\|_{p}\leq N \alpha^{-1/r'}$ so that
$\|m'\chi_{[\alpha,1]}\|_{p}=O(\alpha^{-1/r'})$. \hfill $\Box$
\end{proof}

The following theorem gives sufficient conditions to be satisfied by $m
\in C(0,1]$ to be a multiplier from $AC_{r}$ to $AC_{p},r> p> 1$.

\begin{theor}[\!]
If $r > p>1, m\in L^{v}(I) \cap C(0,1]${\rm ,} where $\frac
{1}{v}=\frac {1}{p}-\frac{1}{r}, m$ is absolutely continuous on
$[\epsilon,1]$ for all $\epsilon>0$ and $\sum _{n=1}^{\infty}(
2^{-n/r'}\|P_{n}m'\|_{p})^{p} <\infty${\rm ,} then $Tg = mg$
defines a multiplier from $AC_{r}$ to $AC_{p}$. Here $P_{n}m' =
m' \cdot \chi_{(2^{- n},2^{-n+1}]}$.
\end{theor}

\begin{proof}
For $g\in AC_{r}, g(0) = 0$. Since $m \in L^{v}(I) \cap C(0,1],
m(t)=O(t^{-1/v})$ in a neighbourhood of 0 so that $|mg(t)| \leq A
\|g'\|_{r}~t^{-1/v+1/r'} = A \|g'\|_{r}t^{1/p'}$, using (1). As $p
> 1, p'<\infty$ we get $\lim_{t\rightarrow 0}mg(t) = 0$ so $mg$
can be continuously extended to $[0,1]$ by assigning $mg(0) = 0$.
Now $(mg)'=m'g+mg'$. Since $m\in L^{v}(I)$ and $g'\in L^{r}(I)$,
it follows that $mg'\in L^{p}(I)$. We show that $m'g\in L^{p}(I)$.
Replacing $g(t)$ by $\int^{t}_{0}g'(s) \d s$ we get for $n =
1,2,3,\dots$,
\begin{align*}
\int^{2^{-n+1}}_{2^{-n}}|m'g|^{p}(t) \d t &\leq
\int^{2^{-n+1}}_{2^{-n}}|m'(t)|^{p} \left(\int ^{t}_{0}|g'(s)| \d s \right)^{p} \d t\\[.2pc]
&\leq \int^{2^{-n+1}}_{2^{-n}}|m'(t)|^{p}(\|g'\|_{r}t^{1/r'})^{p} \d t, \quad \hbox{using (1)},\\[.2pc]
&\leq \|g'\|^{p}_{r} \left(\frac{1}{2^{n-1}}\right)^{p/r'} \int^{2^{-n+1}}_{2^{-n}}|m'(t)|^{p} \d t.
\end{align*}
Adding for $n = 1, 2, 3, \dots$, we get
\begin{equation*}
\|m'g\|_{p}\leq \|g'\|_{r} \left(\sum_{n=1}^{\infty} 2^{-(n-1)p/r'}\int^{2^{-
n + 1}}_{2^{-n}}|m'(t)|^{p} \d t\right)^{1/p} <\infty,
\end{equation*}
since $g'\in L^{r}(I)$ and $\sum_{n=1}^{\infty}(
2^{-n/r'}\|P_{n}m'\|_{p})^{p} < \infty$. Thus $(mg)'\in
L^{p}(I)\subseteq L^{1}(I)$. That $mg$ is absolutely continuous on $I$
follows from Lemma~1 so that $mg \in AC_{p}$. This completes the proof.\hfill $\Box$
\end{proof}

\pagebreak

\begin{remar}$\left.\right.$

{\rm
\begin{enumerate}
\renewcommand\labelenumi{(\arabic{enumi})}
\leftskip .1pc
\item For $f\in L^{r}(I), t^{\frac {1}{r}-\frac
{1}{p}-1}\int^{t}_{0}f(s) \d s \in L^{p}(I)$ for $r\leq p$, see
\cite{6}. The generalized Hardy's inequality is not available in the
case of $r > p$ as in the case of $r = p$, so we get a set of necessary
conditions and another set of sufficient conditions for multipliers in
the case $r > p$.

\item It is easy to see that the condition $\|m'\chi_{[\epsilon,
1]}\|_{p} = O(\epsilon^{-1/r'})$ is equivalent to $\hbox{sup}_{n}2^{-
n/r'}\|P_{n} m'\|_{p} < \infty$.

\item The results of Theorems~3 and 4 compare with those in \cite{2} in view of
Remark~2. In \cite{3}, we get that $m$ is a multiplier from $L^{r}(I)$
to $L^{p}(I), r > p$ if $m\in L^{v}(I)$ and $\sum
_{n=1}^{\infty}(2^{-n/r'}\|P_{n}m'\|_{p})^{v}< \infty$, since
$\frac {1}{v}=\frac {1}{p}-\frac{1}{r}< \frac {1}{p}~,~p < v$ and
$\ell^{p}\subset \ell^{v}$. Thus the results obtained in Theorem~6
are contained in  the results of \cite{3} whereas the necessary
conditions in Theorem~5 compare with those in \cite{3}.
\end{enumerate}

The following examples have bearing on the above said necessary and
sufficient conditions.}
\end{remar}

\setcounter{theore}{0}
\begin{exampl}
{\rm If $m(t)=t^{-1/v+\delta}, 0<\delta<1/v$ then $m(t)\in
L^{v}(I)\cap C(0,1],~m$ is not bounded and $\sum
_{n=1}^{\infty}(2^{-n/r'}\|P_{n}m'\|_{p})^{p}< \infty$. Thus the
sufficiency conditions of Theorem~6 are satisfied so $m$ defines a
multiplier from $AC_{r}$ to $AC_{p},~r>p>1$.}
\end{exampl}

\begin{exampl}
{\rm If $m(t)=t^{-1/v}$ then $m \notin L^{v}(I)$. Also $\sum
_{n=1}^{\infty}(2^{-n/r'}\|P_{n}m'\|_{p})^{p }= \infty$. So $m$ does not
satisfy the sufficiency conditions given by Theorem~6 but
$\int^{1}_{\epsilon}|m'(t)|^{p} \d t = O(\epsilon^{-p/r'})$ so that
necessary condition given by Theorem~5 is satisfied. However, if we take
$g(x)=\int^{x}_{0}t^{-1/r}(1-\ln t)^{-2/r} \d t, x\in I$ then $g\in
AC_{r}$ but $mg \notin AC_{p}$ if $r\leq 2p$ (one can check that
$(mg)'\notin L^{p}(I)$.) So $m$ does not define a multiplier from
$AC_{r}$ to $AC_{p}, p<r\leq 2p$.}
\end{exampl}

\section*{Acknowledgements}

The author is grateful to Prof.~H~L~Vasudeva for many useful discussions
and to Prof.~Ajit Iqbal Singh of University of Delhi for useful
suggestions.


\begin{thebibliography}{9}
\bibitem{1} Baker~J~W, Pym~J~S and Vasudeva~H~L, Totally ordered measure
spaces and their $L^{p}$-algebras, {\it Mathematika} {\bf 29} (1982)
42--54

\bibitem{2} Bhatnagar~S and Vasudeva~H~L, Spaces of multipliers and
their pre-duals for the order multiplication on [0,1], {\it Colloq.
Math.} {\bf 94} (2002) 21--36

\bibitem{3} Bhatnagar~S, Spaces of multipliers and their pre-duals for
the order multiplication on [0,1]~--~II, {\it Colloq. Math.} {\bf 99} (2004)
267--273

\bibitem{4} Larsen~R, An introduction to the theory of multipliers
(Berlin, New York: Springer-Verlag) (1971)

\bibitem{5} Larsen~R, The multipliers of $L^{1}[0,1]$ with order
convolution, {\it Publ. Math. Debrecen} {\bf 23} (1976) 239--248

\bibitem{6} Okikiolu~G~O, Bounded linear transformations in $L^{p}$
space, {\it J. London Math. Soc.} {\bf 41} (1966) 407--414

\bibitem{7} Rickart~C~E, General theory of Banach algebras (New York: D. Van
Nostrand Co.) (1960)

\bibitem{8} Stein~E~M and Weiss~G, Introduction to Fourier analysis on
Euclidean spaces (New Jersey: Princeton) (1971)

\end{thebibliography}
\end{document}